\documentclass[12pt]{amsart}

\usepackage[margin = 1in]{geometry}
\usepackage{amsfonts, amsmath, amscd, amssymb,
latexsym, mathrsfs, mathtools, 
slashed, stmaryrd, verbatim, wasysym, bbm, url,dsfont}
\usepackage[utf8]{inputenc}
\usepackage[backend=bibtex,bibstyle=authoryear,citestyle=alphabetic,giveninits=true,doi=false,isbn=false,url=false,backref]{biblatex}
\usepackage[usenames]{xcolor}
\usepackage{colortbl}
\usepackage{enumerate}
\usepackage[shortlabels]{enumitem}
\usepackage{scalerel}
\usepackage{soul}
\usepackage{diagbox}
\usepackage{array}
\usepackage{nicematrix}
      
\DeclareFieldFormat{labelalphawidth}{\mkbibbrackets{#1}}

\defbibenvironment{bibliography}
  {\list
     {\printtext[labelalphawidth]{%
        \printfield{labelprefix}%
    \printfield{labelalpha}%
        \printfield{extraalpha}}}
     {\setlength{\labelwidth}{\labelalphawidth}%
      \setlength{\leftmargin}{\labelwidth}%
      \setlength{\labelsep}{\biblabelsep}%
      \addtolength{\leftmargin}{\labelsep}%
      \setlength{\itemsep}{\bibitemsep}%
      \setlength{\parsep}{\bibparsep}}%
      }
  {\endlist}
  {\item}
\DeclareNameAlias{author}{last-first}

\newtheorem*{acknowledgements}{Acknowledgements}
\newtheorem*{aistatement}{A.I. usage statement}
\newtheorem{theorem}{Theorem}
\newtheorem{corollary}[theorem]{Corollary}

\newtheorem{lemma}[theorem]{Lemma}

\theoremstyle{definition}

\newtheorem{remark}[theorem]{Remark}

\let\oldsqrt\sqrt
\def\sqrt{\mathpalette\DHLhksqrt}
\def\DHLhksqrt#1#2{%
\setbox0=\hbox{$#1\oldsqrt{#2\,}$}\dimen0=\ht0
\advance\dimen0-0.2\ht0
\setbox2=\hbox{\vrule height\ht0 depth -\dimen0}%
{\box0\lower0.4pt\box2}}

\usepackage{verbatim} %% code from mathabx.sty and mathabx.dcl
\DeclareFontFamily{U}{mathx}{\hyphenchar\font45}
\DeclareFontShape{U}{mathx}{m}{n}{
      <5> <6> <7> <8> <9> <10>
      <10.95> <12> <14.4> <17.28> <20.74> <24.88>
      mathx10
      }{}
\DeclareSymbolFont{mathx}{U}{mathx}{m}{n}
\DeclareFontSubstitution{U}{mathx}{m}{n}
\DeclareMathAccent{\widecheck}{0}{mathx}{"71}

\newcommand\eps\varepsilon
\renewcommand\epsilon\varepsilon

\newcommand{\redd}[1]{}

\newcommand{\abs}[1]{\left\lvert #1 \right\rvert}

\newcommand{\norm}[1]{\lVert #1 \rVert}
\newcommand\inner[1]{\langle #1 \rangle}

\newcommand\sign{\operatorname{sign}}

\newcommand\tr{\operatorname{tr}}

\newcommand\Mand{\text{ and }}

\newcommand\Mwhere{\text{ where }}

\newcommand\paperintro%
        {%
         }
\newcommand\paperbody%
        {%
         }

\newcommand\bbE{\mathbb{E}}

\newcommand\bbL{\mathbb{L}}

\newcommand\bbN{\mathbb{N}}

\newcommand\bbP{\mathbb{P}}

\newcommand\bbR{\mathbb{R}}
\newcommand\bbS{\mathbb{S}}

\newcommand\cH{\mathcal{H}}

\newcommand\cM{\mathcal{M}}

\newcommand\cO{\mathcal{O}}

\DeclareMathAlphabet{\mathpzc}{OT1}{pzc}{m}{it}

\newcommand{\sbs}{\subset}

\DeclareMathOperator{\Span}{span}

\newcommand{\Lpdot}{\bbL^{p(\cdot)}}
\newcommand{\Lpprimedot}{\bbL^{p'(\cdot)}}

\usepackage[refpage]{nomencl}

\makenomenclature

\usepackage{tikz, tikz-cd, tkz-euclide}
\usetikzlibrary{calc,positioning,intersections,quotes,decorations.markings}
\makeatletter
\def\@tocline#1#2#3#4#5#6#7{\relax
  \ifnum #1>\c@tocdepth % then omitth
  \else
    \par \addpenalty\@secpenalty\addvspace{#2}%
    \begingroup \hyphenpenalty\@M
    \@ifempty{#4}{%
      \@tempdima\csname r@tocindent\number#1\endcsname\relax
    }{%
      \@tempdima#4\relax
    }%
    \parindent\z@ \leftskip#3\relax \advance\leftskip\@tempdima\relax
    \rightskip\@pnumwidth plus4em \parfillskip-\@pnumwidth
    #5\leavevmode\hskip-\@tempdima
      \ifcase #1
       \or\or \hskip 1em \or \hskip 2em \else \hskip 3em \fi%
      #6\nobreak\relax
    \hfill\hbox to\@pnumwidth{\@tocpagenum{#7}}\par% <---- \dotfill -> \hfill
    \nobreak
    \endgroup
  \fi}
\makeatother

\def\annu#1{_{% 
  \vbox{\hrule height .2pt 
    \kern 1pt 
    \hbox{$\scriptstyle {#1}\kern 1pt$}% 
  }\kern-.05pt 
  \vrule width .2pt 
}}

\allowdisplaybreaks
\usepackage[colorlinks,citecolor=blue,urlcolor=red,bookmarks=false,hypertexnames=true]{hyperref} 

\makeatletter
\def\keywords{\xdef\@thefnmark{}\@footnotetext}
\makeatother

\title[Boundedness of the Riesz potential associated to diffusion operators]{Boundedness of the Riesz potential associated to diffusion operators in Lebesgue spaces}

\author{Yaozhong Qiu}
\address{UPL, Univ. Paris Nanterre, CNRS, F92000 Nanterre France}
\email{yqiu@parisnanterre.fr}

\author{Cyril Roberto}
\address{UPL, Univ. Paris Nanterre, CNRS, F92000 Nanterre France}
\email{croberto@math.cnrs.fr}

\begin{document}
\begin{abstract}
We revisit boundedness of the Riesz potential associated to a diffusion operator both in fixed and variable exponent Lebesgue spaces. We obtain optimal asymptotics in the former setting and we provide new tools to establish boundedness in the latter. 
\end{abstract}


\maketitle

% Saloff-Coste 
% https://www.scribd.com/document/969148865/Laurent-Saloff-Coste-Aspects-of-Sobolev-Type-Ine-BookZZ-org

% Lieb
% https://sci-hub.fr/10.2307/2007032

\section{Introduction and main results}
In this paper, we prove estimates for the Riesz potential associated to a diffusion operator of the form $L = \Delta - \nabla W \cdot \nabla$ on $\bbR^n$ with $W$ smooth. Denoting by $(P_t)_{t \geq 0}$ the associated semigroup with (infinitesimal) generator $L$, both $L$ and $P_t$ are symmetric by construction in $\bbL^2(\mu)$ where $d\mu = e^{-W}dx$ and $dx$ is Lebesgue measure on $\bbR^n$, i.e. $\int f(Lg)d\mu = \int g(Lf)d\mu$ for all $f, g \in \bbL^2(\mu)$ and similarly for $P_t$. Throughout the paper, $\mu$ is assumed to be a probability measure, i.e. $\int_{\bbR^n} d\mu = 1$. 

The Riesz potential associated to $L$ of parameter $\alpha > 0$ is the operator 
\begin{equation}\label{def:riesz-potential-1}
    I_\alpha = (-L)^{-\alpha}\Pi_0 \ \Mwhere \ \Pi_0f = f - \int_{\bbR^n} fd\mu, 
\end{equation}
that is $I_\alpha$ is the spectral multiplier $\lambda \mapsto (-\lambda)^{-\alpha}$ acting on the subspace of functions with $\mu$-mean zero. Equivalently,
\begin{equation}\label{def:prob-riesz-potential}
    I_\alpha f(x) = \frac{1}{\Gamma(\alpha)} \int_0^\infty t^{\alpha-1}P_t\left(f - \int_{\bbR^n}fd\mu\right)(x)dt, \quad x \in \bbR^n
\end{equation}
where $\Gamma$ is the usual gamma function. The prototypical example of interest considered by this paper is given by the Ornstein-Uhlenbeck operator for which $\smash{W(x) = \frac{\abs{x}^2}{2} - \frac{n}{2}\log(2\pi)}$ and which is symmetric with respect to the $n$-dimensional standard gaussian measure 
\[ d\gamma_n = (2\pi)^{-n/2}e^{-\abs{x}^2/2}dx. \]
It enjoys a spectral decomposition in terms of the Hermite polynomials $(H_k)_{k \in \bbN^n}$ and on each Wiener chaos $\cH_\ell = \Span\{H_k :\abs{k} = \ell\}$ the action of $I_\alpha$ is given by multiplication $I_\alpha \cH_\ell = \ell^{-\alpha}\cH_\ell$ for each $\ell \geq 1$ with the convention $I_\alpha \cH_0 = 0$. Consequently, the gaussian Riesz potential associated to the Ornstein-Uhlenbeck operator can be considered as the gaussian analogue of the euclidean Riesz potential
\begin{equation}\label{def:euc-riesz-potential}
    K_\beta f(x) = c_\beta \int_{\bbR^n} \frac{f(y)}{\abs{y-x}^{n-\beta}}dy, \quad c_\beta = 2^{-\beta}\pi^{-n/2}\frac{\Gamma((n-\beta)/2)}{\Gamma(\beta/2)}  
\end{equation}
defined for $0 < \beta < n$ and $f \in C_0^\infty(\bbR^n)$, since the Fourier transform implies $K_\beta f = (-\Delta)^{-\beta/2}f$, see for instance, \cite[Chapter~V]{stein1970singular}, \cite[\S5.9]{lieb2001analysis}, and \cite[\S1.2]{saloff2002aspects}.

One reason for studying Riesz potentials is that the sharp constant in the $\bbL^p$-$\bbL^q$ estimate 
\[ \norm{K_\beta f}_q \leq C(p, n, \beta) \norm{f}_p, \quad 1 < p < \frac{n}{\alpha} \text{ and } \frac{1}{q} = \frac{1}{p} - \frac{\alpha}{n}, \]
is related to the $\bbL^2$-Sobolev inequality; indeed, the sharp constant (for a specific choice of $p$ and $\beta$) is related to the sharp constant in the Hardy-Littlewood-Sobolev inequality which is known in three cases due to Lieb \cite{lieb1983Sharp}, one of which is dual to the $\bbL^2$-Sobolev inequality. The estimate also implies the $\bbL^p$-Sobolev inequality, $1 < p \neq 2 < n$, through the representation formula 
\[ f(x) = \frac{1}{\omega_{n-1}} \int_{\bbR^n} \frac{\abs{\nabla f(y)}}{\abs{y-x}^{n-1}}dy \] 
where $\omega_{n-1}$ is the volume of the sphere $\bbS^{n-1} \sbs \bbR^n$, see \cite[\S1.2]{saloff2002aspects}. 

For the gaussian Riesz potential, one can only hope for an estimate of the form 
\begin{equation}\label{def:prob-riesz-estimate}
    \norm{I_\alpha f}_{p, \gamma_n} \leq C(n, \alpha, p)\norm{f}_{p, \gamma_n}
\end{equation}
for $1 < p < \infty$. By \cite[Proposition~6.2]{garcia1999spectral} and \cite[\S2]{forzani2000riesz} respectively, it is known $I_\alpha$ is neither of weak type $(1, 1)$ nor is it bounded $\bbL^p(\gamma_n) \rightarrow \bbL^q(\gamma_n)$ with $q > p$, in contrast to the euclidean setting. Meanwhile, $\bbL^p(\gamma_n)$-boundedness is both a consequence of a simple argument exploiting hypercontractivity of the Ornstein-Uhlenbeck semigroup, which we give for completeness below, as well as a special case of an abstract spectral multiplier theorem due to P. A. Meyer \cite[Theorem~3]{meyer1984transformations}, see also \cite[\S6.2]{urbina2019gaussian}. That said, $I_\alpha$ does improve integrability in the same sense that the gaussian measure satisfies a log-Sobolev inequality, see \cite[Theorems~5~\&~6]{bakry1982inegalites1, bakry1982inegalites2}, \cite[Theorem~1]{meyer1982some}, \cite[Theorem~5.7]{karadzhov2005extrapolation}, and \cite[Proposition~8.4]{urbina2019gaussian}. 

% [Meyer-Bakry 1] https://www.numdam.org/item/SPS_1982__16__138_0.pdf
% [Meyer-Bakry 2] https://www.numdam.org/item/SPS_1982__16__146_0.pdf
% [Transformations] https://www.numdam.org/item/SPS_1984__18__179_0.pdf
% [Note] https://www.numdam.org/article/SPS_1982__16__95_0.pdf
% [Some analytic results] https://www.numdam.org/item/RCP25_1982__30__126_0.pdf

Recall the semigroup $(P_t)_{t \geq 0}$ is said to be hypercontractive if for some $\rho > 0$ and all $f \in \bbL^p(\mu)$ it holds 
\begin{equation}\label{def:hypercontractivity}
    \norm{P_tf}_{q, \mu} \leq \norm{f}_{p, \mu} \quad \text{ for all } \quad q \leq 1 + (p - 1)e^{2\rho t}
\end{equation}
for all $1 < p < \infty$. In particular, $(P_t)_{t \geq 0}$ is also a contraction in $\bbL^p(\mu)$, that is $\norm{P_tf}_{p, \mu} \leq \norm{f}_{p, \mu}$. The connection between hypercontractivity and logarithmic Sobolev inequalities was established by Gross \cite{gross1975logarithmic} who proved the Ornstein-Uhlenbeck semigroup is hypercontractive with $\rho = 1$; in contrast, any diffusion semigroup is a contraction in $\bbL^p(\mu)$ for $p \in [1, \infty]$. 

Denoting by $\Pi_\ell$ the projection onto $\cH_\ell$, $\ell \in \bbN_0$, consider for instance $q \geq 2$ and $t_0$ such that $q = 1 + e^{2t_0}$. Then for all $t \geq t_0$,
\begin{align*}
    \norm{P_tf}_{p, \gamma_n}^2 = \norm{P_{t_0}P_{t-t_0}(I - \Pi_0)f}_{p, \gamma_n}^2 \leq \norm{P_{t-t_0}(I - \Pi_0)f}_{2, \gamma_n}^2 &= \sum_{\ell=1}^\infty e^{-2(t - t_0)\ell}\norm{\Pi_\ell f}_{2, \gamma_n}^2 \\
    &\leq e^{-2(t - t_0)}\sum_{\ell=1}^\infty \norm{\Pi_\ell f}_{2, \gamma_n}^2 
\end{align*}
by the semigroup property, hypercontractivity, and orthogonality of the projections respectively. By Parseval's theorem and definition of $t_0$, we find 
\begin{align*}
    \norm{P_tf}_{p, \gamma_n}^2 \leq e^{-2(t - t_0)}\norm{f}_{2, \gamma_n}^2 = (p - 1)e^{-2t}\norm{f}_{2, \gamma_n}^2 \leq (p - 1)e^{-2t}\norm{f}_{p, \gamma_n}^2. 
\end{align*} 
Applying Minkowski's inequality, we obtain 
\begin{align*}
    \norm{I_\alpha f}_{p, \gamma_n} &\leq \frac{1}{\Gamma(\alpha)}\int_0^\infty t^{\alpha-1}\norm{P_tf}_pdt \vphantom{\frac{1}{\Gamma(\alpha)}\left(\int_0^{t_0} t^{\alpha-1}\norm{P_tf}_pdt + \sqrt{p-1}\int_{t_0}^\infty t^{\alpha-1}e^{-t}\norm{f}_pdt\right)} \\
    &\leq \frac{1}{\Gamma(\alpha)}\left(\int_0^{t_0} t^{\alpha-1}\norm{P_tf}_pdt + \sqrt{p-1}\int_{t_0}^\infty t^{\alpha-1}e^{-t}\norm{f}_pdt\right) \\
    &\leq \left(\frac{t_0^\alpha}{\alpha\Gamma(\alpha)} + \sqrt{p-1}\right)\norm{f}_p \vphantom{\frac{1}{\Gamma(\alpha)}\left(\int_0^{t_0} t^{\alpha-1}\norm{P_tf}_pdt + \sqrt{p-1}\int_{t_0}^\infty t^{\alpha-1}e^{-t}\norm{f}_pdt\right)} \\
    &= \left(\frac{\log(p-1)^\alpha}{2^\alpha\Gamma(\alpha+1)} + \sqrt{p-1}\right)\norm{f}_p. \vphantom{\frac{1}{\Gamma(\alpha)}\left(\int_0^{t_0} t^{\alpha-1}\norm{P_tf}_pdt + \sqrt{p-1}\int_{t_0}^\infty t^{\alpha-1}e^{-t}\norm{f}_pdt\right)}
\end{align*}
That is, the above argument leads to the estimate 
\[ C(n, \alpha, p) \leq \frac{\log(p-1)^\alpha}{2^\alpha\Gamma(\alpha+1)} + \sqrt{p-1} \]
which is a bound independent of $n$ and with asymptotic $\sqrt{p}$ as $p \rightarrow \infty$. This feature should not come as a surprise since the gaussian measure is a product measure and the hypercontractivity property is dimension free. Besides the trivial case $p = 2$ in which case $C(n, \alpha, 2) = 1$ is given by the spectral decomposition of the Ornstein-Uhlenbeck operator, we are not aware of any optimal constants in \eqref{def:prob-riesz-estimate}. 

The objectives of this paper are two-fold. The first is to analyse \eqref{def:prob-riesz-estimate} for general diffusion semigroups under the assumption of a (weaker) spectral gap and of hypercontractivity, and to obtain optimal asymptotics, which are apparently new, for the constant $C(n, \alpha, p)$ in the $p \rightarrow \infty$ limit via a modification of the previous argument. The second is to extend the methods developed in the first to the setting of variable exponent spaces, a type of Musielak-Orlicz space. We refer the reader to \S\ref{sec:riesz-boundedness-variable-p} for an introduction to these spaces. This setting is much more involved and the only result available in the literature, to the best of our knowledge, is due to Navas, Pineda, and Urbina \cite[Theorem~1.2]{navas2022boundedness} which deals specifically with the Ornstein-Uhlenbeck operator and exploits the explicit representation of the Mehler heat kernel. Indeed, we show the general strategy alluded above for fixed exponent extends to variable exponents provided one can establish the semigroup is bounded, and we in turn develop some tools to prove such a boundedness property which may also be of independent interest. As application of our methods, we shall generalise the results \cite{teuwen2015note, moreno2019boundedness, navas2022boundedness} to new classes of diffusions and variable exponents. For more references, we refer the reader to \cite{ane2000sur, bakry2004functional, bakry2014analysis} for background on diffusion operators, Poincar\'e inequalities, and hypercontractivity, and the monograph of Urbina \cite{urbina2019gaussian} for discussion of gaussian Riesz potentials. 

The contents of this paper are organised as follows. In \S\ref{sec:riesz-boundedness-fixed-p}, we prove two estimates for the Riesz potential in the fixed exponent setting under the two aforementioned assumptions of a spectral gap and hypercontractivity. In \S\ref{sec:riesz-boundedness-variable-p} we show how to transplant the strategy to the variable exponent setting assuming boundedness of the semigroup; in turn, in \S\ref{sec:semigroup-boundedness-variable-p} we describe how boundedness of the semigroup can be established in the variable exponent setting through pointwise control via the Hardy-Littlewood maximal operator, and finally to complete the picture we provide tools in \S\ref{sec:hlm-boundedness-variable-p} to prove how such pointwise control can be established for some new families of diffusions. 
\section{Boundedness of the Riesz potential in fixed exponent Lebesgue spaces}\label{sec:riesz-boundedness-fixed-p} 

This section aims to prove two results, the first of which is an estimate for the Riesz potential under a spectral gap and the second under hypercontractivity, which are equivalent to a Poincar\'e inequality and a logarithmic Sobolev inequality respectively. Because the former is strictly weaker than the latter, it follows there are examples of operators for which the classical arguments of Meyer based on hypercontractivity do not apply but whose Riesz potentials are nevertheless bounded in $\bbL^p(\mu)$. 

\subsection{Estimates under assumption of a spectral gap}

The key observation which allows us to relax the requirement of hypercontractivity to a spectral gap is that, from the integral definition \eqref{def:prob-riesz-potential} of the Riesz potential, decay of the semigroup in the $\bbL^p$-norm of the following form 
\begin{equation}\label{eq:Lp-exponential-decay}
    \norm{f - \mu(f)}_{p, \mu} \leq K_pe^{-\lambda_pt}\norm{f}_{p, \mu} 
\end{equation}
would immediately provide boundedness of the Riesz potential with an operator norm controlled by $\lambda_p^{-\alpha}$. This, in turn, is a consequence, by a result of Cattiaux, Guillin, and Roberto \cite[Theorem~1.1]{cattiaux2010poincare}, of the $\bbL^2$-Poincar\'e inequality 
\begin{equation}\label{def:poincare-inequality}
    \lambda \norm{f - \mu(f)}_{2, \mu}^2 \leq \norm{\abs{\nabla f}}_{2, \mu}^2
\end{equation}
that we assume holds for any sufficiently smooth $f$ with some constant $\lambda > 0$ independent of $f$. It is of course equivalent to \eqref{eq:Lp-exponential-decay} for $p = 2$, $K_2 = 1$, and $\lambda_2 = \lambda$. Throughout the paper, we denote 
\[ \norm{I_\alpha}_{p \rightarrow p} \vcentcolon= \sup_{\norm{f}_{p, \mu} = 1} \norm{I_\alpha f}_{p, \mu} \]
the operator norm of $I_\alpha$, the dependence on the measure $\mu$ being implicit. 

\begin{theorem}\label{thm:riesz-assuming-poincare}
    Assume $\mu$ satisfies a Poincar\'e inequality \eqref{def:poincare-inequality} with constant $\lambda$. Then for any $p \in (1, \infty)$ and any $\alpha > 0$ it holds 
    \[ \frac{1}{\lambda^\alpha} \leq \norm{I_\alpha}_{p \rightarrow p} \leq \frac{C_p}{\lambda^\alpha} \]
    where $C_p = 2^{1-r_p}/r_p^\alpha$ and $r_p = 2\min(\frac{1}{p}, 1-\frac{1}{p})$. 
\end{theorem}

\begin{proof}
    It turns out that although \cite{cattiaux2010poincare} provides better estimates for $K_p$ and $\lambda_p$, since it would not improve the asymptotic in $p$ of the constant $C_p$ but only the prefactor which is not essential here, we use a weaker form of the results obtained in \cite{cattiaux2010poincare} via interpolation. Indeed, defining $T_t = P_t - \Pi_0$, observe trivially $\norm{T_t}_{1 \rightarrow 1} \leq 2$ and $\norm{T_t}_{\infty \rightarrow \infty} \leq 2$ so that by Riesz-Thorin interpolation with $\norm{T_t}_{2 \rightarrow 2} \leq e^{-\lambda t}$ we immediately obtain \eqref{eq:Lp-exponential-decay} with $K_p = 2^{1-r_p}$ and $\lambda_p = \lambda r_p$ with $r_p$ as defined in the theorem. By Minkowski's inequality, for all $f$ with $\mu$-mean zero, we have 
    \[ \norm{I_\alpha f}_{p, \mu} \leq \frac{1}{\Gamma(\alpha)} \int_0^\infty t^{\alpha-1}\norm{P_tf}_{p, \mu}dt \leq \frac{K_p}{\Gamma(\alpha)} \int_0^\infty t^{\alpha-1}e^{-\lambda_pt}\norm{f}_{p, \mu}dt = \frac{K_p}{\lambda_p^\alpha}\norm{f}_{p, \mu} \]
    which proves the upper bound on $\norm{I_\alpha}_{p \rightarrow p}$. For the lower bound, if $p \in (1, 2]$ then observe if $\lambda$ is optimal in \eqref{def:poincare-inequality} and equality in the Poincar\'e inequality is achieved by an eigenfunction $f \in \bbL^2(\mu)$ associated to the eigenvalue $\lambda$, then $P_tf = e^{-\lambda t}f$ and thus
    \[ I_\alpha f = \frac{1}{\Gamma(\alpha)}\int_0^\infty t^{\alpha-1}P_tfdt = \frac{1}{\Gamma(\alpha)}\int_0^\infty t^{\alpha-1}e^{-\lambda t}fdt = \frac{f}{\lambda^\alpha} \]
    which proves the lower bound after taking $\bbL^p$-norms. In the general case, $\lambda$ still belongs to the spectrum of $-L$ and thus there exists a Weyl sequence of functions $(f_n)$ of unit $\bbL^2(\mu)$-norm satisfying $\norm{(\lambda + L)f_n}_{2, \mu} \rightarrow 0$ as $n \rightarrow \infty$. Since 
    \[ P_tf_n - e^{-\lambda t}f_n = \int_0^t \frac{d}{ds}(e^{-\lambda(t-s)}e^{sL}f_n)ds = \int_0^t e^{-\lambda(t-s)}e^{sL}(\lambda + L)f_nds\]
    and $\norm{P_s}_{p \to p} \leq 1$ for all $p \in [1, \infty]$, taking the 
    $\bbL^p(\mu)$-norm inside the integral by Minkowski's inequality again we obtain 
    \[ 
    \norm{P_tf_n - e^{-\lambda t}f_n}_{p, \mu} 
    \leq 
    \int_0^t e^{-\lambda(t-s)} ds \norm{(\lambda + L)f_n}_{p, \mu} 
    \leq 
    \frac{1}{\lambda} \norm{(\lambda + L)f_n}_{2, \mu} 
    \]
    and thus $\norm{P_tf_n - e^{-\lambda t}f_n}_{p, \mu}$ converges uniformly in $t \geq 0$ to zero as $n \rightarrow \infty$. Thus by the reverse triangle inequality
    \begin{align*}
        \norm{I_\alpha f_n}_{p, \mu} &\geq \frac{1}{\Gamma(\alpha)}\norm{\int_0^\infty t^{\alpha-1}e^{-\lambda t}f_ndt}_{p, \mu} - \frac{1}{\Gamma(\alpha)}\norm{\int_0^\infty t^{\alpha-1}(P_tf_n - e^{-\lambda t}f_n)dt}_{p, \mu} \\
        &\geq \frac{1}{\lambda^\alpha} \norm{f_n}_{p, \mu} - \frac{1}{\Gamma(\alpha)} \int_0^\infty t^{\alpha-1}\norm{P_tf_n - e^{-\lambda t}f_n}_{p, \mu}dt  
    \end{align*}
    which recovers the lower bound in the limit $n \rightarrow \infty$. The case $p \in (2, \infty)$ follows by duality. More precisely, if $f \in \bbL^p(\mu)$, $p > 2$, and $p' = \frac{p}{p-1} \in (1, 2)$ denotes the H\"older conjugate, then 
    \begin{align*}
        \norm{I_\alpha f}_{p, \mu} = \sup_{\norm{g}_{p', \mu} \leq 1} \int (I_\alpha f)gd\mu = \sup_{\norm{g}_{p', \mu} \leq 1} \int f(I_\alpha g)d\mu &\leq \sup_{\norm{g}_{p', \mu} \leq 1} \norm{I_\alpha g}_{p', \mu}\norm{f}_{p, \mu} \\
        &= \norm{I_\alpha}_{p' \rightarrow p'}\norm{f}_{p, \mu} \vphantom{\sup_{\norm{g}_{p', \mu} \leq 1} \int f(I_\alpha g)d\mu}
    \end{align*}
    since $f \in \bbL^p(\mu) \sbs \bbL^2(\mu)$ and we may approximate $g \in \bbL^{p'}(\mu)$ in $\bbL^2(\mu)$. 
    % https://math.stackexchange.com/questions/4759736/infimum-of-numerical-range-infimum-of-spectrum
\end{proof}

\subsection{Estimates under assumption of hypercontractivity}

As mentioned earlier, estimates for the Riesz potential for the gaussian measure follow essentially from hypercontractivity. The point of the following theorem is to prove an improved estimate which is optimal in the gaussian setting (up to constants). 

\begin{theorem}\label{thm:riesz-assuming-hypercontractivity}
    Assume the semigroup $(P_t)_{t \geq 0}$ satisfies the hypercontractivity estimate \eqref{def:hypercontractivity} with constant $\rho > 0$. Then for any $1 < p < \infty$ and any $\alpha > 0$ it holds 
    \[ \norm{I_\alpha}_{p \rightarrow p} \leq \frac{1}{(\rho \epsilon)^\alpha} + \frac{1}{(2\rho(1-\epsilon))^\alpha \Gamma(\alpha+1)} \times \begin{cases} \log(p-1)^\alpha, &p \in [2, \infty) \\ \log(1/(p-1))^\alpha, &p \in (1, 2) \end{cases} \]
    for any $\epsilon \in (0, 1)$. 
\end{theorem}

Before we address the proof, note hypercontractivity \eqref{def:hypercontractivity} implies a Poincar\'e inequality \eqref{def:poincare-inequality} with constant $\lambda = \rho$. 
Thus the above theorem shows that under hypercontractivity one has a bound of the form $\norm{I_\alpha}_{p \rightarrow p} \leq \frac{C_p'}{\lambda^\alpha}$ with $C_p' = \cO((\log p)^\alpha)$, which improves on the asymptotic $C_p \asymp p^\alpha$ obtained under a Poincar\'e inequality. Moreover, as before, no loss occurs at $p = 2$ where the exact estimate $\norm{I_\alpha}_{2 \rightarrow 2} = \frac{1}{\lambda^\alpha}$ is recovered after taking $\rho = \lambda$ and $\epsilon \rightarrow 1$. 

\begin{proof}
    Fix $\epsilon \in (0, 1)$ and $p \geq 2$, the case $p \in (1, 2)$ following again by duality. As mentioned earlier, hypercontractivity implies the Poincar\'e inequality and thus the exponential decay $\norm{P_tf}_{2, \mu} \leq e^{-\rho t}\norm{f}_{2, \mu}$ for all $f$ with $\mu$-mean zero. By Minkowski's inequality, 
    \[ \norm{I_\alpha f}_{p, \mu} \leq \frac{1}{\Gamma(\alpha)} \int_0^{t_0} t^{\alpha-1}\norm{P_tf}_{p, \mu}dt + \frac{1}{\Gamma(\alpha)} \int_{t_0}^\infty t^{\alpha-1}\norm{P_tf}_{p, \mu}dt \]
    where $t_0 = \frac{1}{2\rho(1-\epsilon)}\log(p-1)$. On the first integral, we use that $P_t$ is a contraction so that 
    \[ \int_0^t t^{\alpha-1}\norm{P_tf}_{p, \mu}dt \leq \int_0^T t^{\alpha-1}\norm{f}_{p, \mu}dt = \frac{T^\alpha}{\alpha}\norm{f}_{p, \mu} \]
    while on the second integral we apply hypercontractivity and subsequently the exponential decay property. That is, for $t \geq t_0$ we have $1 + e^{2\rho(1-\epsilon)t} \geq 1 + e^{2\rho(1-\epsilon)t_0} = p$ so that 
    \[ \norm{P_tf}_{p, \mu} = \norm{P_{(1-\epsilon)t}(P_{\epsilon t}f)}_{p, \mu} \leq \norm{P_{\epsilon t}f}_{2, \mu} \leq e^{-\epsilon \rho t}\norm{f}_{2, \rho} \leq e^{-\epsilon \rho t}\norm{f}_{p, \mu}. \]
    Hence 
    \[ \frac{1}{\Gamma(\alpha)}\int_{t_0}^\infty t^{\alpha-1}\norm{P_tf}_{p, \mu}dt \leq \frac{1}{\Gamma(\alpha)}\int_{t_0}^\infty t^{\alpha-1}e^{-\epsilon \rho t}\norm{f}_{p, \mu}dt \leq \frac{1}{(\rho \epsilon)^\alpha}\norm{f}_{p, \mu}. \]
    It follows that 
    \[ \norm{I_\alpha}_{p \rightarrow p} \leq \frac{1}{(\rho \epsilon)^\alpha} + \frac{t_0^\alpha}{\alpha\Gamma(\alpha)} = \frac{1}{(\rho \epsilon)^\alpha} + \frac{\log(p-1)^\alpha}{(2\rho(1-\epsilon))^\alpha \Gamma(\alpha+1)} \]
    as expected. 
    % with $t_0 = \frac{1}{2\rho}\log(p-1)$. On the first integral, we use that $P_t$ is a contraction so that
    % \[ \int_0^t t^{\alpha-1}\norm{P_tf}_{p, \mu}dt \leq \int_0^{t_0} t^{\alpha-1}\norm{f}_{p, \mu}dt = \frac{t_0^\alpha}{\alpha}\norm{f}_{p, \mu} \]
    % while on the second integral we apply hypercontractivity followed by the Poincar\'e contraction property. 
    % Indeed, for $t \geq t_0$ we may write
    % \[ \norm{P_tf}_{p, \mu} = \norm{P_{t_0}P_{t-t_0}f}_{p, \mu} \leq \norm{P_{t-t_0}f}_{2, \mu} \leq e^{-\rho(t - t_0)}\norm{f}_{2, \mu} \leq e^{-\rho(t-t_0)}\norm{f}_{p, \mu} \]
    % and hence 
    % \[ \frac{1}{\Gamma(\alpha)} \int_{t_0}^\infty t^{\alpha-1}\norm{P_tf}_{p, \mu}dt \leq \frac{1}{\Gamma(\alpha)} \int_{t_0}^\infty t^{\alpha-1}e^{-\rho(t-t_0)}\norm{f}_{p, \mu}dt \leq \frac{\sqrt{p-1}}{\rho^\alpha} \norm{f}_{p, \mu}. \]
\end{proof}

The point of the proof is that when we decompose the integral defining the Riesz potential into two parts, there is competition between the contraction and hypercontractivity properties, which give the $p$-dependent term and $p$-independent term in the bound of Theorem \ref{thm:riesz-assuming-hypercontractivity} respectively, the former of which supplies the behaviour of the large and small $p$ asymptotics. Indeed, as we show in the next statement, when we restrict to the Ornstein-Uhlenbeck operator and the standard gaussian measure, the leading term in the $p \rightarrow \infty$ asymptotic comes from the contraction part and therefore corresponds to small $t$. 

\begin{theorem}\label{thm:riesz-optimality-OU}
    Let $\norm{I_\alpha}_{p \rightarrow p}$ be the operator norm of the Riesz potential associated to the Ornstein-Uhlenbeck operator in $\bbR^n$. Then for any $\alpha > 0$ one has 
    \[ \lim_{p \rightarrow 1} \frac{\norm{I_\alpha}_{p \rightarrow p}}{\log(1/(p-1))^\alpha} = \lim_{p \rightarrow \infty} \frac{\norm{I_\alpha}_{p \rightarrow p}}{\log(p -1)^\alpha} = \frac{1}{2^\alpha\Gamma(\alpha+1)}. \]
\end{theorem}

\begin{proof}
    The Ornstein-Uhlenbeck semigroup is hypercontractive with $\rho = 1$ so that by Theorem \ref{thm:riesz-assuming-hypercontractivity} and taking $\epsilon \rightarrow 0$ we have 
    \[ \limsup_{p \rightarrow \infty} \frac{\norm{I_\alpha}_{p \rightarrow p}}{\log(p -1)^\alpha} \leq \frac{1}{2^\alpha\Gamma(\alpha+1)}. \] 
    It follows we need to prove a similar bound for the limit inferior. To that aim, consider the following choice of test function 
    \[ f_0(x) = \mathbbm{1}_{x_1 \geq 0} - \frac{1}{2} = \frac{1}{2}\sign(x_1), \quad x = (x_1, \cdots, x_n) \in \bbR^n. \]
    By the product structure of the Ornstein-Uhlenbeck semigroup and the gaussian measure we observe 
    \[ \norm{I_\alpha f_0}_{p, \gamma_n} = \norm{I_\alpha f_0}_{p, \gamma_1} \Mand \norm{f_0}_{p, \gamma_n} = \norm{f_0}_{p, \gamma_1} = \frac{1}{2} \]
    where with slight abuse of notation we understand $f_0$ both as a function of the ambient variable $x \in \bbR^n$ and of the distinguished variable $x_1 \in \bbR$. With the problem reduced to dimension one, to ease notation we shall take, for the remainder of this proof, all norms with respect to $\gamma_1$ and suppress its appearance in the subscript, i.e.\ $\norm{\cdot}_p = \norm{\cdot}_{p, \gamma_1}$, and interpret $f_0$ as a function of one variable $x \in \bbR$. 

    Starting with 
    \begin{equation}\label{eq:liminf-Ialpha} 
        \liminf_{p \rightarrow \infty} \norm{I_\alpha}_{p \rightarrow p} \geq \liminf_{p \rightarrow \infty} \frac{\norm{I_\alpha f_0}_p}{\norm{f_0}_p} = 2\liminf_{p \rightarrow \infty} \norm{I_\alpha f_0}_p,
    \end{equation}
    the objective is therefore now to bound from below $\norm{I_\alpha f_0}_p$ for $p$ large. Observe that since $f_0$ is odd, $P_tf_0$ is odd so it suffices to consider $x \geq 0$. By Mehler's formula and with $Z$ denoting a one-dimensional standard gaussian random variable, 
    \[ P_tf_0(x) = \bbE(f_0(e^{-t}x + \sqrt{1-e^{-2t}}Z)) = \int_{-\frac{e^{-t}x}{\sqrt{1-e^{-2t}}}}^\infty d\gamma(y) - \frac{1}{2} = \int_0^{\frac{e^{-t}x}{\sqrt{1-e^{-2t}}}} d\gamma(y). \]
    By Fubini's theorem, for any $A > 0$ it holds 
    \begin{align*}
        I_\alpha f_0(x) &= \frac{1}{\Gamma(\alpha)} \int_0^\infty t^{\alpha-1} P_tf_0(x)dt \\ 
        &= \frac{1}{\Gamma(\alpha)} \int_0^\infty \int_0^{\frac{1}{2} \log(1 + x^2/y^2)} dtd\gamma(y) \\
        &= \frac{1}{2^\alpha \Gamma(\alpha+1)} \int_0^\infty \left(\log\left(1 + \frac{x^2}{y^2}\right)\right)^\alpha d\gamma(y) \\ 
        &\geq \frac{1}{2^\alpha \Gamma(\alpha+1)} \left(\log\left(1 + \frac{x^2}{A^2}\right)\right)^\alpha \bbP(0 \leq Z \leq A). 
    \end{align*}
    In turn, since $I_\alpha f_0$ is also odd, 
    \begin{equation}\label{eq:norm-Ialpha}
        \norm{I_\alpha f_0}_p^p = 2\int_0^\infty \abs{I_\alpha f_0}^pd\gamma \geq 2\left(\frac{\bbP(0 \leq Z \leq A)}{2^\alpha \Gamma(\alpha+1)}\right)^p \int_0^\infty \left(\log\left(1 + \frac{x^2}{A^2}\right)\right)^{\alpha p}d\gamma(x)
    \end{equation}
    We study the last integral in detail, using the Laplace method. Let 
    \[ H(x) = \alpha p\log\left(\log\left(1 + \frac{x^2}{A^2}\right)\right) - \frac{x^2}{2}, \quad x \in (0, \infty), \]
    hence
    \begin{align*}
        H'(x) &= \frac{2\alpha p}{A^2} \frac{x}{(1 + \frac{x^2}{A^2})\log(1 + \frac{x^2}{A^2})} - x, \\
        H''(x) &= \frac{2\alpha p}{A^2}\left(\frac{1}{(1 + \frac{x^2}{A^2})\log(1 + \frac{x^2}{A^2})} - \frac{2x^2}{A^2}\frac{1 + \log(1 + \frac{x^2}{A^2}}{((1 + \frac{x^2}{A^2})\log(1 + \frac{x^2}{A^2}))^2}\right) - 1 \\
        &= -\frac{2\alpha p}{A^2} \frac{(-1 + \frac{x^2}{A^2})\log(1 + \frac{x^2}{A^2}) + \frac{2x^2}{A^2}}{((1 + \frac{x^2}{A^2})\log(1 + \frac{x^2}{A^2})^2} - 1
    \end{align*}
    and 
    \[ \int_0^\infty \left(\log\left(1 + \frac{x^2}{A^2}\right)\right)^{\alpha p}d\gamma(x) = \int_0^\infty e^{H(x)} \frac{dx}{\sqrt{2\pi}}. \]
    Observe that $H'(x) = 0$ has a unique solution on the positive axis $x_0$ satisfying $\psi(h(x_0)) = \frac{2\alpha p}{A^2}$, where $\psi(z) = z\log z$ and $h(x) = 1 + \frac{x^2}{A^2}$, and thus $x_0 \asymp \sqrt{2\alpha p/\log p}$ for $p$ large. Fix $\delta > 0$. Expanding to second order, we deduce for any $x \in [x_0 - \delta, x_0 + \delta]$ it holds 
    \[ H(x) \geq H(x_0) - \frac{1}{2}M(x - x_0)^2 \]
    where $M = \max_{y \in [x_0 - \delta, x_0 + \delta]} \abs{H''(y)}$. Hence 
    \begin{align*}
        \int_0^\infty \left(\log\left(1 + \frac{x^2}{A^2} \right)\right)^{\alpha p} d\gamma(x) &\geq 
        \int_{x_0-\delta}^{x_0+\delta} e^{H(x)} \frac{dx}{\sqrt{2\pi}} \\
        &\geq e^{H(x_0)} \int_{x_0 - \delta}^{x_0 + \delta} e^{-\frac{1}{2}M(x - x_0)^2} \frac{dx}{\sqrt{2\pi}} \\
        &= \frac{e^{H(x_0)}}{\sqrt{M}} \int_{-\sqrt{M}\delta}^{\sqrt{M}\delta} e^{-\frac{z^2}{2}} \frac{dz}{\sqrt{2\pi}}.
    \end{align*}
    For large $p$ and thus large $\frac{2\alpha p}{A^2}$, we have $M \rightarrow 2$ near $x_0$ and we may crudely bound $1 \leq M \leq 3$. Returning to \eqref{eq:norm-Ialpha} we have
    \[ \norm{I_\alpha f_0}_p \geq \left(\frac{2}{\sqrt{M}} \int_{-\delta}^\delta e^{-z^2/2}\frac{dz}{\sqrt{2\pi}}\right)^{1/p} e^{H(x_0)/p} \frac{\bbP(0 \leq Z \leq A)}{2^\alpha \Gamma(\alpha + 1)} \]
    and as $p \rightarrow \infty$ using the asymptotic of $x_0 \asymp \sqrt{2\alpha p/\log p}$ we obtain 
    \[ \liminf_{p \rightarrow \infty} \frac{\norm{I_\alpha f_0}_p}{(\log p)^\alpha} \geq \frac{\bbP(0 \leq Z \leq A)}{2^\alpha \Gamma(\alpha + 1)}. \]
    Inserting this into \eqref{eq:liminf-Ialpha} we conclude for any $A > 0$ 
    \[ \liminf_{p \rightarrow \infty} \norm{I_\alpha}_{p \rightarrow p} \geq 2 \liminf_{p \rightarrow \infty} \frac{\norm{I_\alpha f_0}_p}{(\log p)^\alpha} \geq \frac{2\bbP(0 \leq Z \leq A)}{2^\alpha \Gamma(\alpha + 1)} \]
    and taking $A \rightarrow \infty$ yields the desired lower bound. The proof is complete for $p \rightarrow \infty$ while a duality argument gives the result for $p \rightarrow 1$. 
\end{proof}
\section{Boundedness of the Riesz potential in variable exponent Lebesgue spaces}\label{sec:riesz-boundedness-variable-p} 

This section aims to show the previous strategy for fixed exponents, that is combining the contraction property of the semigroup for small time with hypercontractivity and exponential decay in $\bbL^2(\mu)$ for large time, extends to the variable exponent setting. We shall fix a class of variable exponents $p(\cdot): \bbR^n \rightarrow [p_-, p_+]$ measurable with $1 < p_- \leq p_+ < \infty$, satisfying the following two classical regularity conditions.
\begin{enumerate}
    \item \textit{Asymptotic behaviour}: there exists $C_\infty > 0$ and $p_\infty > 1$ such that 
    \[ \abs{p(x) - p_\infty} \leq \frac{C_\infty}{\abs{x}^2}, \quad x \neq 0. \]
    \item \textit{Log-H\"older continuity}: there exists $C_{\log} > 0$ such that for all $x, y \in \bbR^n$ it holds 
    \[ \abs{p(x) - p(y)} \leq \frac{C_{\log}}{\log(e + \abs{x-y}^{-1})}. \]
\end{enumerate}
In the sequel, all variable exponents $p(\cdot)$ are assumed to satisfy these two conditions and are called \emph{regular}. The variable exponent Lebesgue space $\bbL^{p(\cdot)}(\mu)$ associated to $p(\cdot)$ is the set of all measurable functions for which the following norm is finite: 
\[ \norm{f}_{p(\cdot), \mu} \vcentcolon= \inf \left\{\lambda > 0 \,\middle|\, \int \abs{\frac{f(x)}{\lambda}}^{p(x)}d\mu(x) \leq 1\right\}. \]
We refer the reader to the monographs \cite{diening2011lebesgue, cruz2013variable} for more exposition on variable exponent spaces and their properties. 

It should be noted the second condition is called \emph{local} log-H\"older continuity; it usually appears in tandem with log-H\"older continuity \emph{at infinity} which is the first condition but with a denominator of $\log(e + \abs{x})$, and a variable exponent $p(\cdot)$ satisfying both conditions is log-H\"older continuous. In the gaussian setting, despite the blowup near $x = 0$ described in the first condition, such a variable exponent is log-H\"older continuous at infinity and thus in the standard sense. In any case, the motivation for introducing such regularity constraints on the variable exponent apparently came from its sufficiency in establishing the boundedness of the Hardy-Littlewood maximal operator, defined with respect to the Lebesgue measure $dx$, in $\bbL^{p(\cdot)}(dx)$, see for instance the discussion in \cite[\S4]{diening2011lebesgue} and \cite[\S3]{cruz2013variable}. 

For completeness, let us briefly mention that the boundedness $\norm{Tf}_{p(\cdot), \mu} \leq C \norm{f}_{p(\cdot), \mu}$ of an operator $T$ on $\Lpdot(\mu)$ is equivalent to showing 
\[ \int \abs{f(x)}^{p(x)}d\mu \leq 1 \implies \int \abs{Tf(x)}^{p(x)}d\mu \leq C'. \]
Note however that one of the peculiarities of variable exponent spaces is that the boundedness of the operator is not equivalent to the natural expectation that $\int \abs{f(x)}^{p(x)}d\mu \leq C' \int \abs{Tf(x)}^{p(x)}d\mu$ for all $f$; indeed the relationship between the norm and the modular is governed by the norm-modular unit ball property, see \cite[Lemma~3.2.4]{diening2011lebesgue}.  

\begin{theorem}\label{thm:semigroup-boundedness-variable-p-0}
    Let $p(\cdot)$ be a regular variable exponent. Assume the semigroup $(P_t)_{t \geq 0}$ satisfies the hypercontractivity estimate \eqref{def:hypercontractivity} with constant $\rho > 0$ and suppose in addition there exist constants $C_0 \geq 1$, $t_0 > 0$ such that for all $f \in \bbL^{p(\cdot)}(\mu)$ and all $t \in (0, t_0]$ it holds $\norm{P_tf}_{p(\cdot), \mu} \leq C_0\norm{f}_{p(\cdot), \mu}$. Then for any $\alpha > 0$ it holds 
    \[ \norm{I_\alpha f}_{p(\cdot), \mu} \leq \left(\frac{C_0^NT^\alpha}{\Gamma(\alpha+1)} + \frac{2^{3+\alpha-r_{p_-}}}{(\rho r_{p_-})^\alpha}\right)\norm{f}_{p(\cdot), \mu} \]
    where $T = \frac{1}{\rho}\log(\frac{p_+-1}{p_--1})$ and $N = \lceil{T/t_0}\rceil$, i.e.\ $I_\alpha$ is bounded on $\bbL^{p(\cdot)}(\mu)$.
\end{theorem}

An immediate application of this result is the boundedness of the Gaussian Riesz potential, since the $\Lpdot(\mu)$-boundedness of the semigroup for regular variable exponents was established by Moreno, Pineda, and Urbina \cite[Theorem~1.1]{moreno2019boundedness}. 

\begin{proof}
    As before, we pass the $\bbL^{p(\cdot)}(\mu)$-norm inside the integral by Minkowski's inequality and decompose into two parts: 
    \[ \norm{I_\alpha f}_{p(\cdot), \mu} \leq \frac{1}{\Gamma(\alpha)} \int_0^T t^{\alpha-1}\norm{P_tf}_{p(\cdot), \mu} + \frac{1}{\Gamma(\alpha)}\int_T^\infty t^{\alpha-1}\norm{P_tf}_{p(\cdot), \mu}dt \]
    where $T = \frac{1}{\rho}\log(\frac{p_+-1}{p_--1})$ which is twice the time needed for the semigroup to be $\bbL^{p_-}(\mu)$-$\bbL^{p_+}(\mu)$ bounded by hypercontractivity. 
    
    On the first integral, the semigroup property implies for $t \leq T$ and $N = \lceil{T/t_0}\rceil$ that $\norm{P_tf}_{p(\cdot), \mu} \leq C_0^N \norm{f}_{p(\cdot), \mu}$ so that 
    \[ \int_0^T t^{\alpha-1}\norm{P_tf}_{p(\cdot), \mu} \leq C_0^N\norm{f}_{p(\cdot), \mu} = \frac{C_0^NT^\alpha}{\alpha}\norm{f}_{p(\cdot), \mu}. \]
    On the second integral, by \cite[Theorem~3.3.1]{diening2011lebesgue} and hypercontractivity we have 
    \[ \norm{P_tf}_{p(\cdot), \mu} \leq 2\norm{P_tf}_{p_+, \mu} \leq 2\norm{P_{t/2}f}_{p_-, \mu}. \]
    Then, by exponential decay in $\bbL^{p_-}(\mu)$, that is 
    \[ \norm{P_tf}_{p_-, \mu} \leq K_{p_-}e^{-\rho r_{p_-} t}\norm{f}_{p_-, \mu} \]
    with $K_p = 2^{1-r_p}$ and $r_p = 2\min(\frac{1}{p}, 1-\frac{1}{p})$ as given in the proof of Theorem \ref{thm:riesz-assuming-poincare}, we obtain 
    \[ \int_T^\infty t^{\alpha-1}\norm{P_tf}_{p(\cdot), \mu}dt \leq 2^{2-r_p} \norm{f}_{p_-, \mu} \int_T^\infty t^{\alpha-1}e^{-\rho r_{p_-} t/2}dt \leq \frac{2^{2+\alpha-r_{p_-}}}{(\rho r_{p_-})^\alpha} \Gamma(\alpha) \norm{f}_{p_-, \mu} \]
    which gives the theorem since $\norm{f}_{p_-, \mu} \leq 2\norm{f}_{p(\cdot), \mu}$ again by \cite[Theorem~3.3.1]{diening2011lebesgue}. 
\end{proof}

The reader may have noticed that the asymptotic behaviour of the variable exponent $p(\cdot)$ apparently played no role in the proof; indeed, it turns out that, given the hypercontractivity and short-time boundedness of the semigroup in $\Lpdot(\mu)$, the only salient property of $p(\cdot)$ appearing in the proof was that it was bounded between $p_- > 1$ and $p_+ < \infty$. As the next section will show however, the asymptotic decay will play an important role in proving the short-time boundedness. 
\section{Boundedness of the semigroup in variable exponent spaces}\label{sec:semigroup-boundedness-variable-p}

While the literature on hypercontractive semigroups is vast and the examples are many, for instance any diffusion with strictly convex potential $W$, that is $\nabla^2W > 0$ in the sense of quadratic forms, the literature on boundedness of semigroups in variable exponent spaces is, to our best knowledge, currently restricted to the gaussian case studied in \cite{moreno2019boundedness}. Thus in order to give motivation to the previous theorem, this section aims to provide examples of such semigroups and technical tools to prove boundedness. 

Our approach relies on a comparison of the semigroup with the Hardy-Littlewood maximal operator mentioned in the previous section. For a given measure $\mu$, it is defined as 
\[ \cM_\mu f(x) = \sup_{r > 0} \frac{1}{\mu(B(x, r))} \int_{B(x, r)} \abs{f(y)} d\mu(y). \]
While a bound of the form $\abs{P_tf(x)} \leq C M_\mu f(x)$ for all $t > 0$ and $x \in \bbR^n$ is too strong, it turns out that in the gaussian setting it was shown by Teuwen \cite[Theorem~1]{teuwen2015note} it holds for $x$ sufficiently small depending on $t$. The point is that by passing to the Hardy-Littlewood maximal operator we shall demonstrate a general strategy involving decomposition of the space into a compact part in which we may compare $P_tf$ to $\cM_\mu f$ and exploit the boundedness of the latter in $\Lpdot(\mu)$, a result due to Adamowicz, Harjulehto, and H\"ast\"o \cite[Corollary~1.9]{adamowicz2015maximal}, and a noncompact part in which we may use hypercontractivity together with the asymptotic behaviour of $p(\cdot)$. Although we will give a second family of generalisations in the sequel, we start with a warmup to demonstrate the idea of the proof. 

\begin{theorem}\label{thm:semigroup-boundedness-variable-p}
    Let $p(\cdot)$ be a variable exponent satisfying the regularity conditions of \S\ref{sec:riesz-boundedness-variable-p}. Assume that the potential $W$ on $\bbR^n$ satisfies $\rho I \leq \nabla^2W \leq \sigma I$ in the sense of quadratic forms, with $\sigma \geq \rho > 0$. Assume moreover: 
    \begin{enumerate}
        \item there exists $a, C_a > 0$ such that for all $f \in \bbL^1(\mu)$ and all $x \in \bbR^n$ satisfying 
        \[ 1 + \abs{\nabla W(x)} \leq \frac{a}{\sqrt{t}} \]
        it holds $\abs{P_tf(x)} \leq C_a \cM_\mu f(x)$, and 
        \item the various parameters satisfy 
        \[ \frac{B}{a^2} \leq \frac{p_\infty' - 1}{2\rho}, \quad B = \frac{2C_\infty \max(1 + \abs{\nabla W(0)}, \sigma)^2}{(p_--1)(p_\infty-1)}. \]
    \end{enumerate}
    Then there exists a constant $C > 0$ depending on the parameters such that for all $f \in \Lpdot(\mu)$ it holds
    \[ \norm{P_tf}_{p(\cdot), \mu} \leq C \norm{f}_{p(\cdot), \mu}. \]
\end{theorem}

\begin{proof}   
    We introduce the cutoff $\chi_0 = \mathbbm{1}\{1 + \abs{\nabla W(x)} \leq a/\sqrt{t}\}$ and $\chi_1 = 1 - \chi_0$. Moreover, we write $\Phi(f) \lesssim \Phi_2(f)$ to indicate the functional on the left is controlled by the function on the right by a constant independent of $f$ and which depends on the various parameters in the statement of the theorem. 
    
    We first decompose the norm
    \[ \norm{P_tf}_{p(\cdot), \mu} \leq \norm{\chi_0 P_tf}_{p(\cdot), \mu} + \norm{\chi_1 P_tf}_{p(\cdot), \mu} \]
    by the triangle inequality \cite[Theorem~3.2.7]{diening2011lebesgue}, and assume $\int \abs{f}^pd\mu \leq 1$ where the dependence on the spatial variable $x$ has been suppressed in both the function $f$ and the variable exponent $p(\cdot)$ to ease notation. We may also assume without loss of generality $f \geq 0$. For the first term, we have 
    \[ \norm{\chi_0 P_tf}_{p(\cdot), \mu} \lesssim \norm{\chi_0 \cM_\mu f}_{p(\cdot), \mu} \leq \norm{\cM_\mu f}_{p(\cdot), \mu} \lesssim \norm{f}_{p(\cdot), \mu} \]
    since $\mu$ is Radon and $p(\cdot)$ is log-H\"older continuous and so $\cM_\mu$ is bounded on $\Lpdot(\mu)$ according to \cite[Corollary~1.9]{adamowicz2015maximal}. For the second term, let $0 \leq g \in \Lpprimedot(\mu)$ be such that $\norm{g}_{p'(\cdot), \mu} \leq 1$ where $p'(\cdot)$ is pointwise H\"older conjugate to $p(\cdot)$. By the norm-conjugate formula \cite[Corollary~3.2.14]{diening2011lebesgue} we have 
    \begin{align*}
        \norm{\chi_1 P_tf}_{p(\cdot), \mu} \lesssim \sup_{\norm{g}_{p'(\cdot), \mu} \leq 1} \int \abs{\chi_1 P_tf}\abs{g} d\mu &= \sup_{\norm{g}_{p'(\cdot), \mu} \leq 1} \int f P_t(\chi_1 g)d\mu \\ 
        &\lesssim \norm{f}_{p(\cdot), \mu} \sup_{\norm{g}_{p'(\cdot), \mu} \leq 1} \norm{P_t(\chi_1 g)}_{p'(\cdot), \mu}
    \end{align*}
    by approximating $f$ and $g$ with $\bbL^2(\mu)$-functions if necessary, and where the final inequality follows from H\"older's inequality \cite[Lemma~3.2.20]{diening2011lebesgue}. By the triangle inequality again, 
    \[ \norm{P_t(\chi_1 g)}_{p'(\cdot)} \leq \norm{\chi_0 P_t(\chi_1 g)}_{p'(\cdot), \mu} + \norm{\chi_1 P_t(\chi_1 g)}_{p'(\cdot), \mu} \]
    and for $\norm{\chi_0 P_t(\chi_1g)}_{p'(\cdot), \mu}$ we may proceed as before and compare $\chi_0 P_t(\chi_1 g) \lesssim \chi_0 \cM_\mu(\chi_1 g) \leq \cM_\mu g$ and use boundedness of $\cM_\mu$ on $\Lpprimedot(\mu)$ by \cite[Corollary~1.9]{adamowicz2015maximal} since the conjugate exponent $p'(\cdot)$ is also log-H\"older continuous if $p(\cdot)$ is. This implies that 
    \[ \norm{\chi_0 P_t(\chi_1g)}_{p'(\cdot), \mu} \lesssim \norm{\cM_\mu g}_{p'(\cdot), \mu} \lesssim \norm{g}_{p'(\cdot), \mu}. \]
    Finally, to control $\norm{\chi_1 P_t(\chi_1g)}_{p'(\cdot), \mu}$, the point is that $\chi_1$ is supported at infinity where the exponent, according to the asymptotic behaviour assumption, is constant modulo an error linear in $t$ and we can apply hypercontractivity. First note by log-H\"older continuity that the decay condition may be replaced by 
    \[ \abs{p(x) - p_\infty} \leq \frac{C_\infty}{1+\abs{x}^2}, \quad x \in \bbR^n, \]
    and that if $p(\cdot)$ enjoys this asymptotic, so too does $p'(\cdot)$: 
    \[ \abs{p'(x) - p_\infty'} = \abs{\frac{p(x)}{p(x) - 1} - \frac{p_\infty}{p_\infty - 1}} \leq \frac{C_\infty}{(p_- - 1)(p_\infty - 1)(1 + \abs{x}^2)}. \]
    By our assumptions on $W$, observe 
    \[ (1 + \abs{\nabla W(x)})^2 \leq (1 + \abs{\nabla W(0)} + \sigma \abs{x})^2 \leq 2 \max((1 + \abs{\nabla W(0)}), \sigma)^2 (1 + \abs{x}^2) \]
    so that 
    \[ \abs{p'(x) - p_\infty'} \leq \frac{2C_\infty \max((1 + \abs{\nabla W(0)}), \sigma)^2}{(p_- - 1)(p_\infty - 1)(1 + \abs{\nabla W(x)})^2} = \frac{B}{(1 + \abs{\nabla W(x)})^2} \]
    where we set $B$ as in the statement of the theorem. In particular, on the support of $\chi_1$, we have 
    \[ q_1(t) = \max\left(p'_-, p_\infty' - \frac{B}{a^2}t\right) \leq p'(x) \leq \min\left(p'_+, p_\infty' + \frac{B}{a^2}t\right) = q_2(t) \] 
    with $p'_- = (p_+)'$ and $p'_+ = (p_-)'$.

    Assuming $\norm{g}_{p'(\cdot), \mu} \leq 1$, we have $\int \chi_1 \abs{g}^{p'}d\mu \leq \int \abs{g}^{p'}d\mu \leq 1$. Then 
    \begin{align*}
        \int \abs{\chi_1 P_t(\chi_1 g)}^{p'}d\mu \leq \int \chi_1 \abs{\max(1, P_t(\chi_1 g))}^{q_2(t)}d\mu &\leq 1 + \int \abs{P_t(\chi_1 g)}^{q_2(t)}d\mu \\
        &= 1 + \norm{P_t(\chi_1 g)}_{q_2(t), \mu}^{q_2(t), \mu}. \vphantom{\int \chi_1 \abs{\max(1, P_t(\chi_1 g))}^{q_2(t)}}
    \end{align*}
    By hypercontractivity with parameter $\rho > 0$, the assumption $\frac{B}{a^2} \leq \frac{p_\infty' - 1}{2\rho}$ implies 
    \[ \norm{P_t(\chi_1g)}_{q_2(t), \mu}^{q_2(t)} \leq \norm{\chi_1g}_{q_1(t), \mu}^{q_2(t)} \]
    and thus 
    \begin{align*}
        \int \abs{\chi_1 P_t(\chi_1 g)}^{p'}d\mu \leq 1 + \left(\int \abs{\chi_1 g}^{q_1(t)}d\mu\right)^{q_2(t)/q_1(t)} 
        &\leq 1 + \left(\int \abs{\max(1, \chi_1 g)}^{p'}d\mu\right)^{q_2(t)/q_1(t)} \\
        &\leq 1 + \left(1 + \int \abs{g}^{p'}d\mu\right)^{q_2(t)/q_1(t)} \\ 
        &\leq 1 + 2^{q_2(t)/q_1(t)} \vphantom{1 + \left(\int \abs{\chi_1 g}^{q_1(t)}d\mu\right)^{q_2(t)/q_1(t)}}.
    \end{align*}
    It suffices to observe that $(p_+)' \leq q_1(t) \leq q_2(t) \leq (p_-)'$ so that $q_2(t)/q_1(t) \leq (p_-)'/(p_+)'$. Thus we just showed the chain of implications 
    \begin{align*}
        \int \abs{g}^{p'}d\mu \leq 1 \implies \int \abs{\chi_1P_t(\chi_1 g)}^{p'}d\mu \lesssim 1 &\implies \norm{\chi_1 P_t(\chi_1g)}_{p'(\cdot), \mu} \lesssim \norm{g}_{p'(\cdot), \mu} \\ &\implies \norm{\chi_1 P_tf}_{p(\cdot), \mu} \lesssim \norm{f}_{p(\cdot), \mu}. \vphantom{\int \abs{\chi_1P_t(\chi_1 g)}^{p'}d\mu \lesssim 1}
    \end{align*} 
    The theorem follows, i.e. $\norm{P_tf}_{p(\cdot), \mu} \leq C\norm{f}_{p(\cdot), \mu}$ for some $C > 0$. 
\end{proof}
\section{Pointwise semigroup bounds by the Hardy-Littlewood maximal operator}\label{sec:hlm-boundedness-variable-p}

In this final section, our aim is to provide examples of semigroups satisfying a pointwise bound of the form $\abs{P_tf(x)} \lesssim \cM_\mu f(x)$, generically for $x$ depending on $t$. As mentioned earlier, one example in the literature is given by \cite{teuwen2015note} who proved such an estimate in the gaussian setting. We will show in particular that such an estimate persists in the setting of Theorem \ref{thm:semigroup-boundedness-variable-p} as well as with more work in the more general setting of $W = \abs{x}^\alpha$ for $\alpha > 2$; these estimates then supply boundedness of the semigroup in $\bbL^{p(\cdot)}$ and thus of the Riesz potential. We start with two preliminary lemmas. 

Recall that in our setting the semigroup has a smooth kernel $p_t(x, y)$ and that the Hardy-Littlewood maximal operator is defined as 
\[ \cM_\mu f(x) = \sup_{r > 0} \frac{1}{\mu(B(x, r))} \int_{B(x, r)} \abs{f(z)}d\mu(z) \]
where $B(x, r)$ is the ball of radius $r > 0$ centred at $x \in \bbR^n$. For $\theta > 0$, we decompose the space into dyadic annuli
\[ A_k(x, \theta) = \{z \in \bbR^n \mid 2^k\theta \leq \abs{z - x} < 2^{k+1}\theta\} \]
for $k \geq 1$ and $A_0 \vcentcolon= B(x, 2\theta)$. The family $\{A_k(x, \theta)\}_{k \geq 0}$ forms a partition of $\bbR^n$ for all $x \in \bbR^n$ and $\theta > 0$. Define 
\[ a_k(t, x, \theta) = \sup_{z \in A_k} p_t(x, z). \]

\begin{lemma}\label{lem1}
    For all $x \in \bbR^n$ and $\theta > 0$ it holds 
    \[ \abs{P_tf(x)} \leq C \cM_\mu f(x) \]
    with $C = \sum_{k \geq 0} a_k(t, x, \theta)\mu(B(x, 2^{k+1}\theta))$, which may be infinite. 
\end{lemma}

\begin{proof}
    Assume without loss of generality $f \geq 0$ and, to ease notation, write $A_k = A_k(x, \theta)$ and $a_k = a_k(t, x, \theta)$ for fixed $x \in \bbR^n$ and $\theta > 0$. Then 
    \begin{align*}
        P_tf(x) &= \sum_{k \geq 0} \int_{A_k} f(z)p_t(x, z)d\mu(z) \\ 
        &\leq \sum_{k=0}^\infty a_k\int_{A_k} f(z)d\mu(z) \\ 
        &= \sum_{k \geq 0} a_k\mu(B(x, 2^{k+1}\theta)) \left(\frac{1}{\mu(B(x, 2^{k+1}\theta))} \int_{B(x, 2^{k+1}\theta)} f(z)d\mu(z)\right) \\ 
        &\leq \left(\sum_{k=0}^\infty a_k\mu(B(x, 2^{k+1}\theta))\right)\cM_\mu f(x)
    \end{align*}
    as claimed. 
\end{proof}

\begin{lemma}\label{lem2}
    Assume that $W$ is convex. Then for any $x \in \bbR^n$ and any $r > 0$ it holds 
    \begin{equation}\label{eq:volume-ball}
        \mu(B(x, r)) \leq \omega_n r^n e^{-W(x) + r\abs{\nabla W(x)}}
    \end{equation}
    where $\omega_n$ is the euclidean volume of the unit ball in $\bbR^n$. 
\end{lemma}

\begin{proof}
    By convexity, 
    \[ \mu(B(x, r)) = \int_{B(x, r)} e^{-W(y)}dy \leq e^{-W(x)} \int_{B(x, r)} e^{-\inner{\nabla W(x), y-x}}dy \leq e^{-W(x) + r\abs{\nabla W(x)}}\abs{B(x, r)} \]
    which is the desired conclusion since $\abs{B(x, r)} = \omega_n r^n$. 
\end{proof}

\subsection{The bounded hessian case}
Using the previous decomposition, we show how to transfer a heat kernel estimate due to Chen and Wang \cite{chen2024two} to obtain the desired comparison. 

\begin{theorem}\label{thm:example1}
    Assume that $\rho I \leq \nabla^2W \leq \sigma I$ for some $\sigma \geq \rho > 0$. Then for all $t$ sufficiently small and $a > 0$ there exists $C_a > 0$ depending on $n, \rho, \sigma$ such that for all $f \in \bbL^1(\mu)$ and all $x \in \bbR^n$ satisfying
    \[ 1 + \abs{\nabla W(x)} \leq \frac{a}{\sqrt{t}} \]
    it holds 
    \[ \abs{P_tf(x)} \leq C_a \cM_\mu f(x). \]
\end{theorem}

\begin{remark}
    Rewriting the condition as $\sqrt{t} \leq a/(1 + \abs{\nabla W(x)})$, the heuristic is that the diffusion operator $L = \Delta - \nabla W \cdot \nabla$ has two components, one diffusive part (the laplacian) and one drift part (the gradient). In a time interval of length $t$ the diffusive part moves a particle at $x$ by a distance of order $\sqrt{t}$ while the drift will move the particle toward the origin by a distance $\abs{\nabla W(x)}t$. Therefore the drift is not seen by the system if $\abs{\nabla W(x)}t \ll \sqrt{t}$, which is the essence of the bound $\sqrt{t} \leq a/(1 + \abs{\nabla W(x)})$, and thus the system behaves like a pure diffusion with invariant measure the Lebesgue measure. 
\end{remark}

\begin{proof}
    The main ingredient in the proof is a heat kernel estimate for $p_t(x, z)$. Let $L^V = -\frac{1}{2}\Delta + V$ and set $Q_t^V$ for its associated semigroup given by $Q_t^Vf(x) = \int_{\bbR^n} f(z)q_t(x, z)dz$. Set $V(x) = \frac{1}{8}\abs{\nabla W}^2 - \frac{1}{4}\Delta W$, so that the $h$-transform, see \cite[\S1.15.7]{bakry2014analysis}, yields 
    \[ -Lf = -\Delta f + \nabla W \cdot \nabla f = 2e^{W/2}L^V(e^{-W/2}f) \]
    and 
    \[ P_tf(x) = e^{tL}f(x) = e^{W(x)/2}Q_{2t}^V(e^{-W/2}f)(x) = \int_{\bbR^n} f(z)q_{2t}(x, z)e^{\frac{W(x) + W(z)}{2}}d\mu(z), \] 
    i.e. $p_t(x, y) = q_{2t}(x, y)e^{(W(x) + W(y))/2}$ with respect to $\mu$. The assumption $\rho I \leq \nabla^2W \leq \sigma I$ implies by \cite[Theorem~1.1]{chen2024two} that there exist $t_0, c' > 0$ such that for all $x, z \in \bbR^n$ it holds 
    \[ q_t(x, z) \leq c' \frac{e^{-\frac{\abs{x-z}^2}{c't}}}{t^{n/2}} \exp\left(-\frac{1}{c'}\left(t \min(\abs{x}^2, \abs{z}^2) + \abs{x-z}\max(\abs{x}, \abs{z})\right)\right) \leq c' \frac{e^{-\frac{\abs{x-z}^2}{c't}}}{t^{n/2}} \]
    for all $t \in (0, t_0]$. Indeed, $\Delta W = \tr(\nabla^2W) \leq \sigma n$ while $\abs{\nabla W(x)} \geq \rho \abs{x} - \abs{\nabla W(0)}$ so that the Schr\"odinger potential $V$ is bounded below by a quadratic potential (whose heat kernel bounds from above $q_t$). 

    We seek to show the series $\sum_{k \geq 0}^\infty a_k(t, x, \theta) \mu(B(x, 2^{k+1}\theta))$ is convergent. By Lemma \ref{lem2} applied to $\theta = \sqrt{t}$, we know 
    \[ \mu(B(x, 2^{k+1}\sqrt{t})) \leq \omega_n (2^{k+1}\sqrt{t})^n e^{-W(x) + 2^{k+1}\sqrt{t}\abs{\nabla W(x)}} \lesssim t^{n/2}2^{kn}e^{2^{k+1}a}e^{-W(x)} \]
    by the assumption $1 + \abs{\nabla W(x)} \leq a/\sqrt{t}$ and where $\lesssim$ indicates dimension dependent constants which play no role in the $k$-summability. On the other hand, for $t \leq t_0$, we have 
    \[ a_k(t, x, \theta) \leq \sup_{z \in A_k} c' \frac{e^{-\frac{\abs{x-z}^2}{c't}}}{t^{n/2}} \leq c't^{-n/2}e^{-\frac{4^k}{2c'}} \sup_{z \in A_k} e^{\frac{W(x) + W(z)}{2}} \]
    and expanding, 
    \[ W(z) \leq W(x) + \abs{\nabla W(x)}\abs{x - z} + \frac{\sigma}{2} \abs{x - z}^2 \leq W(x) + 2^{k+1}a + 2\sigma 4^kt. \] 
    It follows 
    \[ a_k(t, x, \theta)\mu(B(x, \theta)) \lesssim \exp\left(kn \log 2 + 2^{k+1}a - \frac{4^k}{2c'} + 2^ka + \sigma 4^kt\right) \]
    so that, for $t$ sufficiently small such that $-\frac{4^k}{2c'} + \sigma 4^kt \leq -\frac{4^k}{4c'}$ say, the series is convergent and thus the result follows from Lemma \ref{lem1}. 
\end{proof}

\begin{corollary}
    Assume that $\rho I \leq \nabla^2W \leq \sigma I$ for some $\sigma \geq \rho > 0$, and that $p(\cdot)$ is regular. Then $I_\alpha$ is bounded on $\Lpdot(\mu)$. 
\end{corollary}

\begin{proof}
    By the previous theorem, and by taking $a$ sufficiently large so the second condition of Theorem \ref{thm:semigroup-boundedness-variable-p} is satisfied, we have $\norm{P_tf}_{p(\cdot), \mu} \lesssim \norm{f}_{p(\cdot), \mu}$ for $t$ sufficiently small (thus for $t$ finite and in turn for all $t$ by hypercontractivity) The corollary then follows from Theorem \ref{thm:semigroup-boundedness-variable-p-0}.
\end{proof}

\subsection{The unbounded hessian case}
In this final section, we shall give a preview of how these results can be generalised to a larger family of diffusions. We shall consider the case $W = (1 + \abs{x}^2)^{\alpha/2}$, the regularisation of $\abs{x}^\alpha$ being done for smoothness purposes. In this case, the Schr\"odinger potential $V$ associated to $W$ is bounded below by a constant times $\abs{x}^{2(\alpha-1)} - 1$ so that the function which governs the heat kernel asymptotics in the statement of \cite[Theorem~1.1]{chen2024two} is $g(s) = 1 + s^{2(\alpha-1)}$. 

\begin{theorem}
    Assume that $W = (1 + \abs{x}^2)^{\alpha/2}$ for some $\alpha > 0$. Then for all $t > 0$ and $a > 0$ there exists $C_a > 0$ depending on $n, \alpha$ such that for all $f \in \bbL^1(\mu)$ and all $x \in \bbR^n$ satisfying
    \[ 1 + \abs{\nabla W(x)} \lesssim 1 + \abs{x}^{\alpha-1} \leq \frac{a}{\sqrt{t}} \]
    it holds 
    \[ \abs{P_tf(x)} \leq C_a \cM_\mu f(x). \]
\end{theorem}

\begin{proof}
    Proceeding as before, the measure of the ball with $\theta = 2^{k+1}\sqrt{t}$ is 
    \[ \mu(B(x, 2^{k+1}\sqrt{t})) \lesssim t^{n/2}2^{kn}e^{-W(x) + 2^{k+1}\sqrt{t}\abs{\nabla W(x)}} \lesssim t^{n/2}2^{kn}e^{2^{k+1}a}e^{-W(x)}. \] 
    The heat kernel $q_t$ of the associated Schr\"odinger operator which can be bounded below by $\abs{x}^{2(\alpha-1)} - 1$ modulo constants enjoys by \cite[Theorem~1.1]{chen2024two} the asymptotics
    \[ q_t(x, z) \leq c'\frac{e^{-\frac{\abs{x-z}^2}{c't}}}{t^{n/2}}\exp\left(-\frac{1}{c'}\left(t\min(\abs{x}^{2(\alpha-1)}, \abs{z}^{2(\alpha-1)}) + \abs{x-z}\max(\abs{x}^{\alpha-1}, \abs{z}^{\alpha-1}\right)\right) \] 
    whenever $t \lesssim 1/(1 + \abs{x}^{\alpha-2}) \leq 1/(1 + \min(\abs{x}, \abs{z})^{\alpha-2})$; in particular when $\sqrt{t} \leq a/(1 + \abs{x}^{\alpha-1})$ this condition holds. As before, the contribution of $t^{n/2}$ from the measure of the ball and of $t^{-n/2}$ from the heat kernel cancel, and on $A_k$ we can write 
    \[ \exp\left(-\frac{\abs{x-z}^2}{c't}\right) \leq \exp\left(-\frac{\abs{x-z}^2}{2c't} - \frac{2}{c'}4^k\right). \]
    When $\exp(-\frac{2}{c'}4^k)$ is multiplied by $2^{kn}e^{2^{k+1}a}$, the resulting $k$-summation is finite so finiteness of $\sum_{k \geq 0} a_k(t, x, \theta)\mu(B(x, \theta))$ follows immediately if we show 
    \[ \sup_{z \in A_k} \exp\left(-\frac{1}{c'}\left(\frac{\abs{x-z}^2}{2t} + t\min(\abs{x}^{2\beta}, \abs{z}^{2\beta}) + \abs{x-z}\max(\abs{x}^{\beta}, \abs{z}^{\beta}) + \frac{W(x) - W(z)}{2}\right)\right) \]
    with $\beta = \alpha - 1$ is uniformly bounded above in $t$ and $k$. 

    We can assume $\abs{x} \leq \abs{z}$ otherwise there is nothing to show, since the first three terms in the inner parentheses are nonnegative and the fourth $\frac{1}{2}(W(z) - W(x))$ is nonnegative if $\abs{x} \leq \abs{z}$. To ease notation, we write $\delta = \abs{x-z}$, $r = \abs{x}$, $s = \abs{z}$, so that it suffices to show 
    \[ \inf_{2^k\sqrt{t} \leq \delta \leq 2^{k+1}\sqrt{t}} \frac{\delta^2}{2t} + tr^{2\beta} + \delta s^\beta + \frac{(1 + r^2)^{\alpha/2} - (1 + s^2)^{\alpha/2}}{2} \]
    is uniformly bounded below in $t$. Write $\delta = \lambda s$; we may assume $s \neq 0$ since otherwise there is nothing to show and since $\delta s^\beta = \lambda s^\alpha \gg \frac{1}{2}(1 + s^2)^{\alpha/2}$ for sufficiently large $s$ whenever $\lambda > \frac{1}{2}$, we may also assume $0 \leq \lambda \leq \frac{1}{2} + \epsilon$ for any choice of $\epsilon > 0$, say $\epsilon = \frac{1}{4}$ for concreteness. By the triangle inequality, we have $\delta = \abs{x - z} \geq \abs{z} - \abs{x} = s - r$ so that $\lambda s \geq s - r$ which implies $r \geq (1 - \lambda)s$. Hence the argument of the infimum is bounded below by  
    \begin{equation}\label{lower-bound-1}
        \frac{\lambda^2(1-\lambda)^{2\beta}}{2a^2}s^{2\alpha} + \lambda s^{\alpha} + \frac{(1 + (1-\lambda)^2s^2)^{\alpha/2} - (1 + s^2)^{\alpha/2}}{2} 
    \end{equation}
    since $\sqrt{t} \leq a/(1 + r^\beta)$ and $a^2/t \geq (1 + r^\beta)^2 \geq r^{2\beta} \geq (1 - \lambda)^{2\beta}s^{2\beta}$. Now 
    \begin{align*}
        \frac{(1 + s^2)^{\alpha/2} - (1 + (1 - \lambda)^2s^2)^{\alpha/2}}{2} &\leq \frac{\alpha}{4}(1 + s^2)^{\alpha/2-1}(s^2 - (1 - \lambda)^2s^2) \\
        &\leq \frac{\alpha}{4}(1 - (1 - \lambda)^2)s^2(1 + s^2)^{\alpha/2-1} \\
        &\leq C_\alpha (1 - (1 - \lambda)^2)(1 + s^\alpha) \vphantom{\frac{\alpha}{2}}
    \end{align*}
    for some $C_\alpha > 0$ and so \eqref{lower-bound-1} is bounded below by a quadratic in $v = s^\alpha$ of the form $c_2v^2 + c_1v + c_0$ where
    \[ c_2 = \frac{\lambda^2(1 - \lambda)^{2\beta}}{2a^2}, \quad c_1 = \lambda - C_\alpha(1 - (1 - \lambda)^2), \quad c_0 = -C_\alpha(1 - (1 - \lambda)^2). \]
    This quadratic achieves a minimum value of $-c_1^2/(4c_2) + c_0$ and it remains to uniformly bound this as a function of $\lambda \in [0, \frac{3}{4}]$. This is a continuous function away from zero so it suffices to check its behaviour near zero where 
    \[ -\frac{c_1^2}{4c_2} + c_0 = \frac{\lambda^2(1 - 2C_\alpha)^2}{2\lambda^2/a^2} - C_\alpha + \cO(\lambda^2) = \frac{a^2(1 - 2C_\alpha)^2}{2} - C_\alpha + \cO(\lambda^2) \]
    which implies $-c_1^2/(4c_2) + c_0$ is continuous on all of $\lambda \in [0, \frac{3}{4}]$ and thus bounded. This implies that $\sum_{k \geq 0} a_k(t, x, \theta)\mu(B(x, \theta)) < \infty$ and the theorem follows. 
\end{proof}

\begin{corollary}
    Assume that $W = (1 + \abs{x}^2)^{\alpha/2}$ for some $\alpha > 2$. Assume $p(\cdot)$ is regular and in addition 
    \[ \abs{p(x) - p_\infty} \leq \frac{C_\infty}{(1 + \abs{x}^{\alpha-1})^2}, \quad x \neq 0. \]
    Then $I_\alpha$ is bounded on $\Lpdot(\mu)$. 
\end{corollary}

\begin{proof}
    Note $(P_t)_{t \geq 0}$ is hypercontractive (in fact ultracontractive), see for instance \cite[\S7.7]{bakry2014analysis}. By the previous theorem, we have a comparison between the semigroup and the Hardy-Littlewood maximal operator for the region in which $1 + \abs{x}^{\alpha-1} \leq a/\sqrt{t}$. We may follow the proof of Theorem \ref{thm:semigroup-boundedness-variable-p} mutatis mutandis to conclude the semigroup is $\bbL^{p(\cdot)}(\mu)$-bounded; that $p(\cdot)$ enjoys the assumed decay implies $p'(\cdot)$ as well modulo constants. On the other hand, in the part of the argument where we use hypercontractivity, that is the complement $1 + \abs{x}^{\alpha-1} > a/\sqrt{t}$, the additional decay assumption on $p'(\cdot)$ ensures $p'(\cdot)$ is again at most a linear error away from $p_\infty'$ which can be made sufficiently small by taking $a$ sufficiently large. Boundedness of the Riesz potential in $\bbL^{p(\cdot)}(\mu)$ in turn follows from Theorem \ref{thm:semigroup-boundedness-variable-p-0}. 
\end{proof}

\begin{acknowledgements}
\textup{This work has been conducted within the FP2M federation (CNRS FR 2036) and FSM CY Initiative. This project has received funding from the European Union’s Horizon 2020 research and innovation programme under the Marie Sk\l{}odowska-Curie grant agreement No 101034255. \scalerel*{\includegraphics{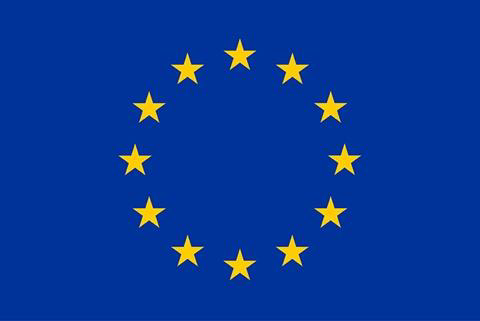}}{A}}
\end{acknowledgements}

\begin{aistatement}
\textup{The authors declare artificial intelligence was not used to develop the mathematical content nor writing of this paper.} 
\end{aistatement}

\printbibliography

\end{document}